\documentclass[12pt]{article}
\usepackage{amssymb}

\usepackage{sw20lart}

\input{tcilatex}
\begin{document}

\author{Alan Horwitz}
\title{Compositions of Polynomials with Coefficients in a given Field}
\date{8/20/01}
\maketitle

\begin{abstract}
Let $F\subset K$ be fields of characteristic $0,$ and let $K[x]$ denote the
ring of polynomials with coefficients in $K$. Let $p(x)=\sum%
\limits_{k=0}^{n}a_{k}x^{k}\in K[x],a_{n}\neq 0.$ For $p\in K[x]\backslash
F[x],$ define $D_{F}(p),$ the $F$ deficit of $p$, to equal $n-\max \{0\leq
k\leq n:a_{k}\notin F\}.$ For $p\in F[x],$ define $D_{F}(p)=n$. Let $%
p(x)=\sum\limits_{k=0}^{n}a_{k}x^{k}$, $q(x)=\sum%
\limits_{j=0}^{m}b_{j}x^{j}, $ with $a_{n}\neq 0,$ $b_{m}\neq 0$, $%
a_{n},b_{m}\in F$,\textit{\ }$b_{j}\notin F$ for some $j\geq 1.$ Suppose
that $p\in K[x]$, $q\in K[x]\backslash F[x],\;p$ not constant. Our main
result is that $p\circ q\notin F[x]$ and $D_{F}(p\circ q)=D_{F}(q)$. With
only the assumption that $a_{n}b_{m}\in F,$ we prove the inequality $%
D_{F}(p\circ q)\geq D_{F}(q).$ This inequality also holds if $F$ and $K$ are
only rings. Similar results are proven for fields of finite characteristic
with the additional assumption that the characteristic of the field does not
divide the degree of $p$. Finally we extend our results to polynomials in
two variables and compositions of the form $p(q(x,y)),$ where $p$ is a
polynomial in one variable.
\end{abstract}

\section{Introduction}

\footnotetext{
Key words: polynomial, field, composition, iterate}Suppose that $p$ and $q$
are polynomials such that their composition, $p\circ q,$ has all rational
coefficients. Must the coefficients of $p$ or $q$ be all rational ? The idea
for this paper actually started with the following more general question.
Let $F\subset K$ be fields of characteristic $0,$ and let $K[x]$ denote the
ring of polynomials with coefficients in $K$. Suppose that $p$ and $q$ are
polynomials in $K[x],$ and $p\circ q\in F[x].$ Must $p$ or $q$ be in $F[x]$
? The answer is yes (see Theorem \ref{T3}) if the leading coefficient and
the constant term of $q$ are each in $F$. Theorem \ref{T3} follows easily
from a more general result (Theorem \ref{T1}) concerning the $F$ deficit,
denoted by $D_{F},$ of the composition of two polynomials. $D_{F}$ is
defined as follows: If $p\in K[x]\backslash F[x]$, $\deg (p)=n$, let $x^{r}$
be the largest power of $x$ with a coefficient \textbf{not} in $F$. We
define the $F$ deficit of $p,$ $D_{F}(p)$, to be $n-r$. For $p\in F[x],$
define $D_{F}(p)=n$. For example, if $F=Q$ (rational numbers), $K=$ $R$
(real numbers), and $p(x)=x^{5}-5x^{3}+\sqrt{3}x^{2}-x+1$, then $D_{F}(p)=3$%
. Now suppose that the leading coefficients of $p$ and $q$ are in $F$, and
that some coefficient of $q$(other than the constant term) is \textbf{not}
in F(so that $q\notin F[x]$). Our main result, Theorem \ref{T1}, states that 
$D_{F}(p\circ q)=D_{F}(q)$. With the weaker assumption that only the \textit{%
product }of the leading coefficients of $p$ and $q$ is in $F$ we prove the 
\textit{inequality} $D_{F}(p\circ q)\geq D_{F}(q)$ (see Theorem \ref{T2}).
It is interesting to note that if $q\in F[x],$ then we get the different
equality $D_{F}(p\circ q)=D_{F}(p)D_{F}(q)$.

We also prove (Theorem \ref{T4}) some results about the deficit of the 
\textbf{iterates}, $p^{[r]},$ of $p$ which require less assumptions than
those of Theorem \ref{T1}. In particular, $D_{F}(p^{[r]})=D_{F}(p)$ with
only the assumption that the leading coefficient of $p$ is in $F$. This
assumption is necessary in general as the example $p(x)=ix$ shows with $F=R$
and $K=C$ (complex numbers).

One can, of course, define the $F$ deficit for \textit{any} two sets $%
F\subset K.$ While Theorem \ref{T1} does not hold in general if $F$ and $K$
are not fields, we can again prove the inequality $D_{F}(p\circ q)\geq
D_{F}(q)$ if $F$ and $K$ are rings (see Theorem \ref{ring}).

For fields of \textbf{finite characteristic }$d$, Theorem \ref{T1} follows
under the additional assumption that $d$ does not divide $\deg (p)$.

Finally we extend our results to polynomials in two variables(using a
natural definition of $D_{F}$ in that case) and compositions of the form $%
p(q(x,y)),$ where $p$ is a polynomial in one variable. Our proof easily
extends to compositions of the form $p(q(x_{1},...,x_{r})).$ However, the
analog of Theorem \ref{T1} does not hold in general for compositions of the
form $p(q_{1}(x,y),q_{2}(x,y))$ (even when $q_{1}=q_{2}$), where $p$ is also
a polynomial in two variables.

There are connections between some of the results in this paper and earlier
work of the author in \cite{H} and \cite{HR}, where we asked questions such
as: If the composition of two power series, $f$ and $g$, is even, must $f$
or $g$ be even ? One connection with this paper lies in the following fact:
If $F=$ $R$ and $K=$ $C$, then $F[x]$ is invariant under the linear operator 
$L(f)(z)=\bar{f}(\bar{z})$. Of course, the even functions are invariant
under the linear operator $L(f)(z)=f(-z)$. Note that in each case $L\circ
L=I $. This connection does not extend to fields $F$ in general, however,
since such an operator $L$ may not exist. The methods and results we use in
this paper are somewhat similar to those of \cite{H} and \cite{HR}, but
there are some key differences. Also, we only consider polynomials in this
paper, since there is really no useful notion of the $F$ deficit for power
series which are not polynomials.

\textit{Acknowledgment}: We thank the referee for suggesting the notion of
the $F$ deficit and its use in strengthening the original version of Theorem 
\ref{T1}.

\section{Main Results}

Let $F\subset K$ be sets, with $F[x]$ equal to the set of all polynomials
with coefficients in $F$.

\begin{definition}
Let $p(x)=\sum\limits_{k=0}^{n}a_{k}x^{k}\in K[x],a_{n}\neq 0.$ For $p\in
K[x]\backslash F[x],$ define $D_{F}(p),$ the $F$ deficit of $p$, to equal $%
n-\max \{0\leq k\leq n:a_{k}\notin F\}.$ For $p\in F[x],$ define $D_{F}(p)=n$%
.
\end{definition}

Note that $D_{F}(p)=n$ if and only if $a_{k}\in F$ $\forall k\geq 1$ and $%
D_{F}(p)=0$ if and only if $a_{n}\notin F$.

Most of the results in this paper concern the case when $F$ and $K$ are 
\textbf{fields}.

We shall need the following easily proven properties. For any fields $%
F\subset K$

\begin{equation}
u\in F,\;v\in K\backslash F\Rightarrow uv\in K\backslash F(\text{if }u\neq 0)%
\text{ and }u+v\in K\backslash F  \label{e0}
\end{equation}

and for fields of characteristic $0$

\begin{equation}
v\in K\backslash F\Rightarrow nv\in K\backslash F\text{ for any }n\in Z_{+}
\label{e1}
\end{equation}

Assume for the rest of this section that $F$ is a proper nonempty \textbf{%
subfield} of $K$, which is a field\textbf{\ }of \textbf{characteristic zero}%
. Later in the paper we discuss the case where $K$ is a field of finite
characteristic or just a ring.

The following result shows that, under suitable assumptions, $q$ and $p\circ
q$ have the same $F$ deficit.

\begin{theorem}
\label{T1}Suppose that $p(x)=\sum\limits_{k=0}^{n}a_{k}x^{k}\in K[x]$, $p$
not constant, $q(x)=\sum\limits_{j=0}^{m}b_{j}x^{j}\in K[x]\backslash F[x]$
with $a_{n}\neq 0,$ $b_{m}\neq 0$, $a_{n},b_{m}\in F$,\textit{\ }$%
b_{j}\notin F$ for some $j\geq 1.$ Then $p\circ q\notin F[x]$ and $%
D_{F}(p\circ q)=D_{F}(q)$.
\end{theorem}

\proof%
Let $d=D_{F}(q)<m.$ Since $b_{m}\in F,$ $d\geq 1.$ By the definition of $%
D_{F}$, $b_{m-d}\notin F,$ but $b_{m-(d-1)},...,b_{m}\in F$. Also, since $p$
is not constant, $n\geq 1.$ We have

\begin{equation}
(p\circ q)(x)=\sum\limits_{k=0}^{n}a_{k}\left(
\sum\limits_{j=0}^{m}b_{j}x^{j}\right) ^{k}  \label{comp}
\end{equation}

Consider the coefficient of $x^{mn-d}$ in $(p\circ q)(x).$ Since $%
mn-d>mn-m=m(n-1),$ this coefficient will only arise from the summand above
with $k=n,$ namely $a_{n}(q(x))^{n},$ which equals

\begin{equation}
a_{n}\left( \sum_{\;i_{0}+\cdots +i_{m}=n}\dfrac{n!}{(i_{0})!\cdots (i_{m})!}%
(b_{0})^{i_{0}}\cdots (b_{m}x^{m})^{i_{m}}\right)  \label{1}
\end{equation}

To get an exponent of $mn-d$ in (\ref{1}), $\sum_{k=0}^{m}ki_{k}=mn-d$.
Along with $\sum_{k=0}^{m}i_{k}=n$ this implies

\begin{equation}
mi_{0}+(m-1)i_{1}+\cdots +i_{m-1}=d  \label{2}
\end{equation}
Note that since $b_{j}\notin F$ for some $j\geq 1,$ $d<m$, which implies
that $m-(d+1)\geq 0.$ Now $mi_{0}+(m-1)i_{1}+\cdots +(d+1)i_{m-(d+1)}>d$ if
some $i_{j}\neq 0$ for $0\leq j\leq m-(d+1).$ That proves

\begin{equation}
i_{j}=0\text{ for }0\leq j\leq m-(d+1)  \label{3}
\end{equation}

By (\ref{2}) and (\ref{3}), $di_{m-d}+(d-1)i_{m-(d-1)}+\cdots +i_{m-1}=d$.
Since $b_{j}\in F$ for $j\geq m-(d-1),$ the only way to get a coefficient in
(\ref{1}) \textit{not} in $F$ is if $i_{m-d}\neq 0,$ which implies that $%
i_{m-d}=1,$ $i_{m-d+1}=i_{m-d+2}=\cdots =i_{m-1}=0.$ Also, from $i_{m-d}+$ $%
i_{m-d+1}+\cdots +i_{m}=n$ we have $i_{m}=n-1$. Hence the only way to obtain 
$x^{mn-d}$ in (\ref{1}) using $b_{m-d}$ is $n\left( b_{m-d}x^{m-d}\right)
^{1}\left( b_{m}x^{m}\right) ^{n-1}.$ Now $b_{m-d}b_{m}^{n-1}\notin F$(by (%
\ref{e0})), and all of the other terms in (\ref{1}) which contribute to the
coefficient of $x^{mn-d}$ involve $b_{m-(d-1)},...,b_{m}.$ Hence, by (\ref
{e0}) and (\ref{e1}), the coefficient of $x^{mn-d}$ in (\ref{1}) is\textbf{\
not} in $F$, and it follows that $p\circ q\notin F[x]$. Now we want to show
that the coefficient of $x^{r}$ in (\ref{comp}) will lie in $F$ if $r>mn-d.$
Write $r=mn-d^{\prime },$ where $d^{\prime }<d.$ Since $mn-d^{\prime }>mn-d,$
this coefficient will only arise in (\ref{comp}) with $k=n.$ Arguing as
above, to get an exponent of $mn-d^{\prime }$ in (\ref{1}), it follows that $%
i_{j}=0$ for $0\leq j\leq m-(d^{\prime }+1)$. Since $m-(d^{\prime }+1)\geq
m-d,$ $i_{j}=0$ for $0\leq j\leq m-d,$ which implies that the coefficient of 
$x^{r}$ in (\ref{1}) only involves $b_{k}$ with $k>m-d$. Since $%
b_{m-(d-1)},...,b_{m}\in F,$ the coefficient of $x^{r}$ in (\ref{comp}) is
also in $F$, and thus $D_{F}(p\circ q)=d$. 
\endproof%

If $q\in K[x]/F[x]$ and $b_{0}\in F,$ then $b_{j}\notin F$ for some $j\geq
1. $ Theorem \ref{1} then implies

\begin{corollary}
\label{C1}Suppose that $p(x)=\sum\limits_{k=0}^{n}a_{k}x^{k}\in K[x]$, $p$
not constant, $q(x)=\sum\limits_{j=0}^{m}b_{j}x^{j}\in K[x]\backslash F[x]$.
Suppose that $a_{n}\neq 0,$ $b_{m}\neq 0$, $a_{n},b_{m},b_{0}\in F.$ Then $%
p\circ q\notin F[x]$ and $D_{F}(p\circ q)=D_{F}(q)$.
\end{corollary}

\begin{example}
Let $F=Q,K=R,$ $p(x)=x^{3}+2x^{2}-\sqrt{2}x+1$,

$q(x)=x^{2}+\sqrt{3}x+5$. Then

\begin{eqnarray*}
p(q(x)) &=&\allowbreak x^{6}+3\sqrt{3}x^{5}+26x^{4}+37\sqrt{3}x^{3}+\left( -%
\sqrt{2}+146\right) x^{2}+ \\
&&\allowbreak \left( 95\sqrt{3}-\sqrt{2}\sqrt{3}\right) x+176-5\sqrt{2}
\end{eqnarray*}
. Hence $D_{F}(p\circ q)=1=D_{F}(q)$.
\end{example}

Theorem \ref{T1} assumes that $q\notin F[x].$ For $q\in F[x]$ we have

\begin{theorem}
\label{T1a}Suppose that $p(x)=\sum\limits_{k=0}^{n}a_{k}x^{k}\in K[x]$, $%
q(x)=\sum\limits_{j=0}^{m}b_{j}x^{j}\in F[x],$ with $a_{n}\neq 0,$ $%
b_{m}\neq 0.$ Then $D_{F}(p\circ q)=D_{F}(p)D_{F}(q)$.
\end{theorem}

\proof%
If $p$ is constant, then $p\circ q$ is constant, and thus $D_{F}(p\circ
q)=0=D_{F}(p)D_{F}(q)$. So assume now that$\;p$ is not constant.

\textit{Case 1:} $a_{n}\in F$ and $p\notin F[x]$

Let $d=D_{F}(p)\Rightarrow d>0,a_{n-d}$ $\notin F,$ and $a_{n-d+1},...,a_{n}%
\in F.$ Consider the coefficient of $x^{mn-md}$ in $(p\circ q)(x).$ This
coefficient will only arise in (\ref{comp}) with $k\geq n-d$. Since $%
a_{k}\in F$ for $k>n-d,$ the only way to get a coefficient not in $F$ is
with $a_{n-d}(q(x))^{n-d}=a_{n-d}b_{m}^{n-d}x^{mn-md}+\cdots $. Since $%
a_{n-d}b_{m}^{n-d}\notin F$, the coefficient of $x^{mn-md}$ is not in $F$.
It also follows easily that if $r>mn-md,$ then the coefficient of $x^{r}$ in
(\ref{comp}) is in $F$. Thus $D_{F}(p\circ q)=mn-(mn-md)=md=D_{F}(p)D_{F}(q)$%
.

\textit{Case 2:} $a_{n}\notin F$

Then $D_{F}(p)=0$ and $a_{n}b_{m}^{n}\notin F\Rightarrow D_{F}(p\circ
q)=0=D_{F}(p)D_{F}(q)$.

\textit{Case 3:} $p\in F[x]$

Then $D_{F}(p\circ q)=mn=D_{F}(p)D_{F}(q)$.

\endproof%

\begin{example}
Let $F=Q,K=R,$ $p(x)=x^{4}-\sqrt{2}x$, and $q(x)=x^{2}+3x$. Then $%
p(q(x))=x^{8}+12x^{7}+54x^{6}+108x^{5}+81x^{4}-\sqrt{2}x^{2}-3\sqrt{2}x$.
Hence $D_{F}(p\circ q)=6=(3)(2)=D_{F}(p)D_{F}(q).$
\end{example}

\begin{remark}
Theorem \ref{T1a} implies that if $q\in F[x],$ then $D_{F}(p\circ q)\geq
D_{F}(q)$.
\end{remark}

\begin{remark}
Theorem \ref{T1} does not hold in general if $a_{n}$ and/or $b_{m}$ are not
in $F$. For example, let $p(x)=\sqrt{2}x^{3}+x^{2}-x+\sqrt{5}$, $q(x)=3\sqrt{%
2}x^{2}+\sqrt{3}x+5$, where $F=Q,K=R$. Then $D_{F}(p\circ q)=1$and $%
D_{F}(q)=0$, and thus $D_{F}(p\circ q)\neq D_{F}(q).$ However, with the
weaker assumption that $a_{n}b_{m}\in F,$ one can prove an \textit{inequality%
}.
\end{remark}

\begin{theorem}
\label{T2}Suppose that $p(x)=\sum\limits_{k=0}^{n}a_{k}x^{k}\in K[x]$, $p$
not constant, $q(x)=\sum\limits_{j=0}^{m}b_{j}x^{j}\in K[x],$ with $%
a_{n}\neq 0,$ $b_{m}\neq 0$, $a_{n}b_{m}\in F.$ Then $D_{F}(p\circ q)\geq
D_{F}(q)$.
\end{theorem}

\proof%
. \textit{Case 1: }$q\notin F[x]$ and $b_{m}\in F$.

By (\ref{e0}), $a_{n}\in F$ as well. If $b_{j}\notin F$ for some $j\geq 1,$
then by Theorem \ref{T1}, $D_{F}(p\circ q)=D_{F}(q)$. Now suppose that $%
b_{j}\in F$ for all $j\geq 1$. It is not hard to show that the coefficient
of any power of $x>m(n-1)$ cannot involve $b_{0}$, and hence $D_{F}(p\circ
q)\geq mn-m(n-1)=m=D_{F}(q)$.

\textit{Case 2: }$q\notin F[x]$ and $b_{m}\notin F$. Then $D_{F}(q)=0$ and
the inequality follows immediately.

\textit{Case 3: }$q\in F[x].$ Then $D_{F}(p\circ q)\geq D_{F}(q)$ by Theorem 
\ref{T1a}(see the remark following the proof). 
\endproof%

\begin{remark}
Theorem \ref{T2} does not hold in general if $a_{n}b_{m}\notin F.$ For
example, let $F=Q,K=R,$ $p(x)=\sqrt{2}x^{3}+x^{2}-x+1$, and $q(x)=x^{2}+%
\sqrt{3}x+5.$ Then clearly $D_{F}(p\circ q)=0$ while $D_{F}(q)=1$.
\end{remark}

As an application of Theorem \ref{T1} we have the following result. Note
that we do \textbf{not} assume that the leading coefficient of $p$ is in $F$.

\begin{proposition}
\label{P1}Suppose that $p(x)=\sum\limits_{k=0}^{n}a_{k}x^{k}\in K[x]$, $p$
not constant, $q(x)=\sum\limits_{j=0}^{m}b_{j}x^{j}\in K[x]\backslash F[x],$
with $a_{n}\neq 0,\;b_{m}\neq 0,$ and $b_{m}\in F.\;$If $b_{j}\notin F$ for
some $j\geq 1,$ then $p\circ q\notin F[x].$
\end{proposition}

\proof%
If $a_{n}\notin F,$ then $a_{n}b_{m}^{n}\notin F$, which implies that $%
p\circ q\notin F[x]$ since $a_{n}b_{m}^{n}$ is the coefficient of $x^{mn}$
in $p\circ q$. If $a_{n}\in F$, then $p\circ q\notin F[x]$ by Theorem \ref
{T1}. 
\endproof%

\begin{lemma}
\label{L1}Suppose that $q(x)=\sum\limits_{j=0}^{m}b_{j}x^{j}\in F[x]$ and $%
p\circ q\in F[x]$, $p(x)=\sum\limits_{k=0}^{n}a_{k}x^{k},a_{n}\neq
0,\;b_{m}\neq 0,q$ not constant. Then $p\in F[x].$
\end{lemma}

\proof%
Note that $D_{F}(q)=m\geq 1>0.$ Then by Theorem \ref{T1a}, $D_{F}(p)=\dfrac{%
D_{F}(p\circ q)}{D_{F}(q)}=\dfrac{mn}{m}=n$, and thus $a_{k}\in F$ for $%
k\geq 1$. Since $p\circ q\in F[x],$ $p(q(0))=\sum%
\limits_{k=0}^{n}a_{k}b_{0}^{k}\in F$. Since $b_{0}\in F$, this implies that 
$a_{0}\in F$. Hence $p\in F[x].$ 
\endproof%

Now we answer the following question mentioned in the introduction. Suppose
that $p\circ q\in F[x].$ Must $p$ or $q$ be in $F[x]$ ?

\begin{theorem}
\label{T3}Suppose that $p,q\in K[x]$ with $p\circ q\in F[x],$ $%
q(x)=\sum\limits_{j=0}^{m}b_{j}x^{j},a_{n}\neq 0,b_{m}\neq 0,b_{0},b_{m}\in
F.$ Then $p\in F[x]$ or $q\in F[x]$. In addition, if $p\circ q$ is not
constant, then $p\in F[x]$ \textit{and }$q\in F[x]$.
\end{theorem}

\proof%
Suppose $p\circ q\in F[x].$ If $p\circ q$ is constant$,$ then $p$ and/or $q$
is constant. If $p(x)=c$, then $(p\circ q)(x)=c$, which implies that $c\in F$
and hence $p\in F[x]$. If $q(x)=c$, then $c\in F$ since $b_{0}\in F$ and
hence $q\in F[x]$. Now suppose that $p\circ q$ is not constant. Then $q$ is
not constant. If $q\notin F[x]$, then $b_{j}\notin F$ for some $j\geq 1.$ By
Proposition \ref{P1}, $p\circ q\notin F[x]$, a contradiction. Hence $q\in
F[x]$. Lemma \ref{L1} then shows that $p\in F[x]$ as well. 
\endproof%

\begin{remark}
Note that no restriction is needed on the leading coefficient of $p.$
However, some restriction on the \textbf{leading coefficient} and \textbf{%
constant} term of $q$ are needed in order for Theorem \ref{T3} to hold in
general. Simple examples are $p(x)=x-c,q(x)=x+c$ or $p(x)=\dfrac{1}{c}%
x,q(x)=cx$, with $c\in K,$ $c\notin F$.
\end{remark}

\begin{remark}
Theorem \ref{T3} does not hold in general if $F$ equals the \textit{%
complement} of a field. For example, if $F=$ irrational numbers, let $%
p(x)=x^{2}$, $q(x)=\pi x^{2}+x+\pi $. Then neither $p$ nor $q$ has all
irrational coefficients, and the leading coefficient and constant term of $q$
are irrational. However, $p(q(x))=\allowbreak \pi ^{2}x^{4}+2\pi
x^{3}+\left( 2\pi ^{2}+1\right) x^{2}+2\pi x+\pi ^{2}$, which has all
irrational coefficients.
\end{remark}

\begin{remark}
If $S$ is any subset of $K$(not necessarily a subfield), we say that $S$ is
a \textbf{deficit set} if Theorem \ref{T1} holds with $F$ replaced by $S$
throughout. For example, if $K=C=$ complex numbers, then it is not hard to
show that $S=R\cup I=$ set of all real or imaginary numbers is a deficit
set. It would be interesting to determine exactly what a deficit set must
look like for a given field $K$.
\end{remark}

\subsection{Iterates}

We now prove the analogs of Theorems \ref{T1}, \ref{T2}, and \ref{T3} when $%
p=q$. In this case we require less assumptions. In particular, for the
analog of Theorem \ref{T2}, we require no assumptions whatsoever. Let $%
p^{[r]}$ denote the $r$th iterate of $p$.

\begin{theorem}
\label{T4}Suppose that $p(x)=\sum\limits_{k=0}^{n}a_{k}x^{k}\in
K[x]\backslash F[x]$, with $a_{n}\neq 0$,$\;a_{n}\in F.$ Then, for any
positive integer $r,$ $p^{[r]}\notin F[x]$ and $D_{F}(p^{[r]})=D_{F}(p)$.
\end{theorem}

\proof%
Note that if $n=0$, then $a_{0}\in F\Rightarrow p\in F[x].$ Hence $n\geq 1.$
If $n=1$, then $p(x)=a_{1}x+a_{0},$\ $a_{1}\in F,a_{0}\notin F.$ It is not
hard to show that 
\[
p^{[r]}(x)=(a_{1})^{r}x+a_{0}\sum_{k=0}^{r-1}(a_{1})^{k}\; 
\]

Now $a_{0}\sum\limits_{k=0}^{r-1}(a_{1})^{k}\notin F$ since $%
\sum\limits_{k=0}^{r-1}(a_{1})^{k}\in F$. Hence $p^{[r]}(x)\notin F[x]$ and $%
D_{F}(p^{[r]})=1=D_{F}(p).$ Assume now that $n\geq 2$. First we prove the
theorem for $p\circ p,$%
\begin{equation}
(p\circ p)(x)=\sum_{k=0}^{n}a_{k}\left( \sum_{j=0}^{n}a_{j}x^{j}\right) ^{k}
\label{compit}
\end{equation}

If $a_{j}\notin F$ for some $j\geq 1,$ then $D_{F}(p\circ p)=D_{F}(p)$ by
Theorem \ref{T1} with $p=q$. So suppose now that $a_{j}\in F$ for $j\geq 1$
and $a_{0}\notin F.$ First let $k=n$ in (\ref{compit}) to get

\begin{equation}
a_{n}\left( \sum_{\;i_{0}+\cdots +i_{n}=n}\dfrac{n!}{(i_{0})!\cdots (i_{n})!}%
(a_{0})^{i_{0}}\cdots (a_{n}x^{n})^{i_{n}}\right)  \label{4}
\end{equation}

It follows easily that the highest power of $x$ in (\ref{4}) involving $%
a_{0} $ is $n(n-1),$ obtained by letting $i_{0}=1,i_{j}=0$ for $2\leq j\leq
n-1,$ $i_{n}=n-1.$ The coefficient of $x^{n(n-1)}$ in (\ref{4}) is $%
na_{0}a_{n}^{n-1}\notin F$ by (\ref{e0}) and (\ref{e1}). The only other way
to obtain $x^{n(n-1)}$ is by letting $k=n-1$ in (\ref{compit}) and letting $%
i_{n}=n-1$ in $a_{n-1}\left( \sum\limits_{\;i_{0}+\cdots +i_{n}=n-1}\dfrac{%
(n-1)!}{(i_{0})!\cdots (i_{n})!}(a_{0})^{i_{0}}\cdots
(a_{n}x^{n})^{i_{n}}\right) $. This gives a coefficient of $x^{n(n-1)}$
which does \textit{not} involve $a_{0}$. Hence the coefficient of $%
x^{n(n-1)} $ in $p\circ p$ equals $na_{0}a_{n}^{n-1}+c,$ where $c\in F.$ By (%
\ref{e0}), $na_{0}a_{n}^{n-1}+c\notin F.$ Finally, it is not hard to show
that any power of $x$ in (\ref{compit}) greater than $n(n-1)$ cannot involve 
$a_{0}.$ Thus $D_{F}(p\circ p)=n^{2}-n(n-1)=n=D_{F}(p).$ Now consider $%
p^{[r]}=p^{[r-2]}\circ q$ where $r\geq 3,$ and $q=p\circ
p=\sum\limits_{j=0}^{m}b_{j}x^{j},$ $m=n^{2}$. Since $D_{F}(p\circ
p)=D_{F}(p)\leq n$, $D_{F}(p\circ p)<n^{2}$ since $n\geq 2$. Hence $%
b_{j}\notin F$ for some $j\geq 1.$ Since $b_{m}=a_{n}^{n+1}\in F$ and the
leading coefficient of $p^{[r-2]}$ is also in $F,$ $%
D_{F}(p^{[r]})=D_{F}(p^{[r-2]}\circ q)=D_{F}(q)$(by Theorem \ref{T1}$%
)=D_{F}(p\circ p)=D_{F}(p).$ It also follows that $p^{[r]}\notin F[x]$ since 
$D_{F}(p^{[r]})=D_{F}(p)\leq n<n^{2}.$ 
\endproof%

\begin{remark}
If $f(x)=\dfrac{x}{ax-1},$ then $f(f(x))=x\in F(x)=$ ring of formal power
series in $x$. However, $f\notin F(x)$ if $a\notin F,$ which implies that
the first part of Theorem \ref{T4} fails in general for formal power
series(we have not defined $D_{F}(f)$ for $f\in F(x)$).
\end{remark}

\begin{remark}
Theorem \ref{T4} is not simply a trivial application of Theorem \ref{T1}
using induction on $r,$ with $q=p^{[r-1]}$. The reason is that one requires $%
b_{j}\notin F$ for some $j\geq 1$ to apply Theorem \ref{T1}.
\end{remark}

\begin{example}
Let $F=R,K=C,$ and $p(x)=x^{3}+4x^{2}-3ix+2i$. Then $p(p(x))=\allowbreak
x^{9}+12x^{8}+\left( 48-9i\right) x^{7}+\left( 68-66i\right) x^{6}+\left(
5-96i\right) x^{5}+\allowbreak \left( -8+72i\right) x^{4}+\left(
132-56i\right) x^{3}+\left( -84-2i\right) x^{2}+\allowbreak \left(
39+36i\right) x-10-6i\Rightarrow D_{F}(p\circ p)=2=D_{F}(p).$
\end{example}

We now prove an inequality which holds for \textbf{all} $p$ in $K[x]$.

\begin{theorem}
Let $p\in K[x].$ Then $D_{F}(p^{[r]})\geq D_{F}(p)$.
\end{theorem}

\proof%
. If $p\in F[x],$ then $p^{[r]}\in F[x],$ which implies that $%
D_{F}(p^{[r]})=n^{r}\geq n=D_{F}(p)$. If $p\in K[x]\backslash F[x]$ and $%
a_{n}\in F,$ then by Theorem \ref{T4}, $D_{F}(p^{[r]})=D_{F}(p)$. Finally,
if $a_{n}\notin F,$ then $D_{F}(p)=0\leq D_{F}(p^{[r]}).$ 
\endproof%

We now prove the analog of Theorem \ref{T3} for iterates.

\begin{theorem}
\label{T5}Suppose that $p\in K[x],$ $p(x)=\sum%
\limits_{k=0}^{n}a_{k}x^{k},a_{n}\neq 0,\;a_{n}\in F.$ Assume also that $%
p^{[r]}\in F[x]$ for some positive integer $r.$ Then $p\in F[x]$.
\end{theorem}

\proof%
If $p\notin F[x],$ then $p^{[r]}\notin F[x]$ by Theorem \ref{T4}. 
\endproof%

\begin{remark}
Theorem \ref{T5} does not hold in general if $a_{n}\notin F.$ For a
counterexample, if there exists $a\in F$ with $a^{1/r}\notin F$, then let $%
p(x)=a^{1/r}x$.\footnote{%
If $F=$ algebraic numbers and $K=$ real numbers, then such an $a$ does not
exist.}
\end{remark}

\section{Several Variables}

As earlier, assume throughout that $F$ is a proper nonempty subfield of $K$,
which is a field\textbf{\ }of characteristic zero. We now extend the
definition of the $F$ deficit to polynomials in two variables. Write $%
p(x,y)=\sum\limits_{k=0}^{n}p_{k}(x,y),$ where each $p_{k}$ is homogeneous
of degree $k,p_{n}\neq 0$. If $p\in K[x,y]\backslash F[x,y],$ define $%
D_{F}(p)=n-\max \{k:p_{k}\notin F[x,y]\}.$ For $p\in F[x,y],$ define $%
D_{F}(p)=n$. Then Theorems \ref{T1} and \ref{T2}, with similar assumptions,
do \textbf{not} hold in general for compositions of the form $p(q(x),q(x))$,
where $q$ is a polynomial in one variable and $p$ is a polynomial in two
variables. For example, let $F=R,$ $K=C$, $p(x,y)=x^{2}-y^{2}+1$, $%
q(x)=x^{2}+ix$. Then $p(q(x),q(x))=\allowbreak 1$ and thus $%
D_{F}(p(q,q))=0<1=D_{F}(q).$ Indeed, Theorems \ref{T1} and \ref{T2} even
fail for iterates of the form $p(p(x,y),p(x,y)).$ For example, let $F=Q,$ $%
K= $ $R$, and $p(x,y)=y^{2}-x^{2}+\sqrt{3}x-\sqrt{5}y$. Then $%
p(p(x,y),p(x,y))=\sqrt{3}y^{2}-\sqrt{3}x^{2}+3x-\sqrt{3}\sqrt{5}y-\sqrt{5}%
y^{2}+\sqrt{5}x^{2}-\sqrt{5}\sqrt{3}x+\allowbreak 5y,$ which implies that $%
D_{F}(p(p,p))=0<1=D_{F}(p).$

However, we can prove similar theorems for compositions of the form $%
p(q(x,y)),$ where $p$ is a polynomial in \textbf{one} variable.

\begin{theorem}
\label{2var}Suppose that $p(x)=\sum\limits_{k=0}^{n}a_{k}x^{k}\in K[x],0\neq
a_{n}\in F,$ $p$ not constant. Suppose that $q\in K[x,y]\backslash F[x,y]$, $%
q(x,y)=\sum\limits_{j=0}^{m}q_{j}(x,y),$ where each $q_{j}$ is homogeneous
of degree $j$ with $0\neq q_{m}\in F[x,y]$. If $q_{j}(x,y)\notin F[x,y]$ for
some $j\geq 1$, then $p\circ q=p(q(x,y))\notin F[x,y]$ and $D_{F}(p\circ
q)=D_{F}(q)$.
\end{theorem}

\proof%
Our proof is very similar to the proof of Theorem \ref{T1}, except that we
have to work with the homogeneous polynomials $q_{j}(x,y)$ instead of the
monomials $x^{j}.$ This only complicates things a little. 
\begin{equation}
p(q(x,y))=\dsum\limits_{k=0}^{n}a_{k}\left(
\dsum\limits_{j=0}^{m}q_{j}(x,y)\right) ^{k}  \label{co2}
\end{equation}

Let $d=D_{F}(q)<m.$ By the definition of $D_{F}$, $q_{m-d}\notin F[x,y],$

$q_{m-(d-1)},...,q_{m}\in F[x,y]$. Also, $p$ not constant $\Rightarrow n\geq
1$ and $q_{n}\in F[x,y]\Rightarrow d>0.$ Now $(q_{j}(x,y))^{k}$ is
homogeneous of degree $jk$, and $k<n$ implies that $jk\leq j(n-1)\leq
m(n-1)<mn-d$. Hence a term of degree $mn-d$ can only arise in (\ref{co2}) if 
$k=n$, which gives $a_{n}(q(x,y))^{n}=\left(
\sum\limits_{j=0}^{m}q_{j}(x,y)\right) ^{n}=$

\begin{equation}
a_{n}\left( \dsum\limits_{\;i_{0}+\cdots +i_{m}=n}\dfrac{n!}{(i_{0})!\cdots
(i_{m})!}(q_{0})^{i_{0}}\cdots (q_{m})^{i_{m}}\right)  \label{10}
\end{equation}

Note that $m-d\geq 1\Rightarrow m-(d+1)\geq 0.$ Arguing exactly as in the
proof of Theorem \ref{T1}, to get an exponent of $mn-d$ in (\ref{10})

\begin{equation}
i_{j}=0\text{ for }0\leq j\leq m-(d+1)  \label{11}
\end{equation}

Thus the only way to get a coefficient in (\ref{10}) \textit{not} in $F$ is
if $i_{m-d}\neq 0,$ which implies that $i_{m-d}=1,$ $i_{m-d+1}=i_{m-d+2}=%
\cdots =i_{m-1}=0.$ Also, from $i_{m-d}+$ $i_{m-d+1}+\cdots +i_{m}=n$ we
have $i_{m}=n-1$. (\ref{10}) then becomes $na_{n}q_{m-d}q_{m}^{n-1}$, which
we shall now show has at least one coefficient not in $F.$ Let $%
g=q_{m-d}q_{m}^{n-1},$ which is homogeneous of degree $mn-d$. Write $%
q_{m}^{n-1}(x,y)=\sum\limits_{k=0}^{m(n-1)}c_{k}x^{k}y^{m(n-1)-k},$ $%
q_{m-d}(x,y)=\sum\limits_{r=0}^{m-d}b_{r}x^{r}y^{m-d-r}.$ Note that $%
c_{k}\in F$ for all $k$, while $b_{r}\notin F$ for some $r$. Let 
\[
M=\max \{r:0\leq r\leq m-d,b_{r}\notin F\},N=\max \{k:0\leq k\leq
m(n-1),c_{k}\neq 0\} 
\]
Clearly $M$ and $N$ are well defined, $b_{M}\notin F$, and $c_{N}\in F.$
Consider the coefficient of $x^{M+N}y^{mn-d-M-N}$ in $g.$ One way to obtain
this coefficient is $\left( b_{M}x^{M}y^{m-d-M}\right) \left(
c_{N}x^{N}y^{m(n-1)-N}\right) =b_{M}c_{N}x^{M+N}y^{mn-d-M-N}$ $.$ There are
other ways to obtain this coefficient if $N>0$ and $M<m-d.$ Since $c_{k}=0$
for $k>N$, one must choose $c_{k}x^{k}y^{m(n-1)-k}$ from $q_{m}^{n-1}$ with $%
k<N$ and $b_{r}x^{r}y^{m-d-r}$ from $q_{m-d}$ with $r>M$, which all involve
coefficients in $F.$ Since $b_{M}c_{N}\notin F$, the coefficient of $%
x^{M+N}y^{mn-d-M-N}$ in $g$ is not in $F.$ Thus the coefficient of $%
x^{M+N}y^{mn-d-M-N}$ in $na_{n}q_{m-d}q_{m}^{n-1}$ is not in $F$, which
implies that $.p(q(x,y))\notin F[x,y]$.

Now write $(p\circ q)(x,y)=\sum\limits_{l=0}^{mn}h_{l}(x,y),$ where each $%
h_{l}$ is homogeneous of degree $l$. Again, arguing exactly as in the proof
of Theorem \ref{T1}, since $q_{m-(d-1)},...,q_{m}\in F[x,y],$ it follows
that $h_{l}\in F[x,y]$ for $l>mn-d.$ This implies that $D_{F}(p\circ
q)=mn-d. $ 
\endproof%

\begin{remark}
Theorem \ref{2var} can be easily extended to compositions of the form $%
p(q(x_{1},...,x_{r})).$
\end{remark}

\section{Rings}

Theorem \ref{T1} does not hold in general if $F$ is just\textit{\ }a \textit{%
ring}. For example, if $F=Z$, the ring of integers and $K=Q$, let $%
p(x)=x^{2}+\dfrac{2}{3}x$ and $q(x)=6x^{2}+\dfrac{3}{2}x$. Then $%
a_{2},b_{2}, $ and $b_{0}$ are in $Z$, and $p(q(x))=\allowbreak \allowbreak
36x^{4}+18x^{3}+\dfrac{25}{4}x^{2}+x$, which implies that $2=D_{F}(p\circ
q)\neq D_{F}(q)=1$. Theorem \ref{T2} also does not hold if $F$ is a ring.
However, if $F\subset K,$ where $F$ and $K$ are rings of finite or infinite
characteristic, we can prove

\begin{theorem}
\label{ring}Suppose that $p(x)=\sum\limits_{k=0}^{n}a_{k}x^{k}\in K[x]$, $p$
not constant,$\;q(x)=\sum\limits_{j=0}^{m}b_{j}x^{j}\in K[x]\backslash F[x],$
with $a_{n}\neq 0,$ $b_{m}\neq 0$, $a_{n},b_{m}\in F$,\textit{\ }$%
b_{j}\notin F$ for some $j\geq 1.$ Then $D_{F}(p\circ q)\geq D_{F}(q)$.
\end{theorem}

\proof%
. Letting $d=D_{F}(q)$, the proof follows exactly as in the proof of Theorem 
\ref{T1}, except that we cannot conclude that $b_{m-d}b_{m}\notin F$ if $F$
is only a ring. However, it does still follow that the coefficient of $x^{r}$
in (\ref{comp}) will lie in $F$ if $r>mn-d.$ Hence, even if $b_{m-d}b_{m}\in
F,$ it follows that $D_{F}(p\circ q)\geq d$. 
\endproof%

\section{ Fields of Finite Characteristic}

Theorem \ref{T1} also does not hold in general if the field $F$ has \textit{%
finite }characteristic. For example, suppose that $K$ is a finite field of
order $4$, $F=Z_{2}\subset K$. Let $p(x)=x^{2}$, $q(x)=x^{2}+3x.$ Then $%
p(q(x))=\allowbreak x^{4}+(3+3)x^{3}+(3\times 3)x^{2}=\allowbreak
x^{4}+2x^{2}$. Thus $\allowbreak D_{F}(q)=1$ while $D_{F}(p\circ q)=2$ . The
problem here is that the characteristic of $K$ divides the degree of $p$. If
we assume that this does not happen, we have

\begin{theorem}
\label{FF}Let $F\subset K$ be fields of characteristic $t$. Suppose that $%
p(x)=\sum\limits_{k=0}^{n}a_{k}x^{k}\in K[x]$, $p$ not constant, $%
q(x)=\sum\limits_{j=0}^{m}b_{j}x^{j}\in K[x]\backslash F[x],$ with $%
a_{n}\neq 0,$ $b_{m}\neq 0$, $a_{n},b_{m}\in F$,\textit{\ }$b_{j}\notin F$
for some $j\geq 1.$ If $t\nmid n$, then$\;p\circ q\notin F[x]$ and $%
D_{F}(p\circ q)=D_{F}(q)$.
\end{theorem}

\proof%
We need the fact that if $r\in Z_{+}$ with $r<t,$ then $ru\neq 0$ for any $%
u\in K.$ This easily implies that $nu\neq 0$ if $t\nmid n$. It follows that
if $u\notin F$, then $nu\notin F$ if $t\nmid n$. Hence, letting $d=D_{F}(q)$
and $u=b_{m-d}b_{m}\notin F$ we have $nb_{m-d}b_{m}\notin F$. Now the proof
follows exactly as in the proof of Theorem \ref{T1}$.$ 
\endproof%

One can also prove versions of Theorems \ref{T2} and \ref{T3} for fields of
finite characteristic. Theorem \ref{T3} also requires the additional
assumption that $t\nmid n.$ 
\endproof%

\section{Applications}

The main theorems in this paper give information about the coefficients of $%
p\circ q,$ and about the coefficients of the iterates of $p$. All of the
examples we give here use $F=$ rationals, $K=$ reals, though of course it is
possible to construct examples from other fields of characteristic $0,$ from
finite fields, or from rings. For example, let $p(x)=x^{2}+c,$ where $c$ is
irrational. By Theorem \ref{T4}, $D_{F}(p)=2\Rightarrow D_{F}(p^{[r]})=2$
for any $r\in Z_{+}$, which implies that the coefficient of $x^{2^{r}-2}$ in 
$p^{[r]}(x)$ is irrational, while the coefficient of $x^{2^{r}-1}$ in $%
p^{[r]}(x)$ must be rational.

Also, suppose that, given $r(x)\in K[x],$ one wants to determine if
nonlinear polynomials $p,q\in K[x]$ exist with $r=p\circ q$. Given $p$ or $q$
as well, Theorems \ref{T1} or \ref{T3} can sometimes be used to give a quick
negative answer. For example, let $r(x)=x^{6}+ax^{5}+bx^{4}+\cdots ,$ where $%
a$ is rational and $b$ is irrational, and $q(x)=x^{3}+Bx^{2}+\cdots $, where 
$B$ is irrational$.$ If $r=p\circ q$, then the leading coefficient of $p$
equals $1,$ and by Theorem \ref{T1}, $D_{F}(q)=D_{F}(r)=2$. But $D_{F}(q)=1$
and thus no such $p$ exists.

The applications given here are probably of limited value. It would be nice
to find other, perhaps more useful, applications of the theorems in this
paper.

\section{Entire Functions}

The obvious extension of $F[x]$ to the class of\textit{\ entire} functions $%
E $ is 
\[
S_{F}=\{f\in E:f(z)=\sum_{k=0}^{\infty }a_{k}z^{k},\;a_{k}\in F\;\forall k\} 
\]

While there is no reasonable notion of $D_{F}(f)$ when $f$ is not a
polynomial, one can attempt to extend Theorem \ref{T3} to $E$. The question
then becomes: Suppose that $f(z)$ is entire and $q(z)$ is a polynomial, with
leading coefficient and constant term in $F$. If $f\circ q\in S_{F},$ must $%
f\in S_{F}$ or $q\in S_{F}$ ? The following theorem gives a negative answer
to this question for a large class of fields $F$.

\begin{theorem}
\label{Tr2}Let $F$ be a subfield of $C,$ with either $F=R$ or $\pi
^{2}\notin F$. Then there exists an entire function $f(z)$ and a polynomial $%
q(z)=a_{2}z^{2}+a_{1}z+a_{0}$ such that :

(1) $f\notin S_{F}$ and $q\notin S_{F}$

(2) $a_{0}$ and $a_{2}$ are both in $F$

(3) $f\circ q\in S_{F}$
\end{theorem}

\proof%
Case 1: $F=$ $R$

Let $f(z)=\cos (i\pi \sqrt{z+2i})=\cosh (\pi \sqrt{z+2i})$ and $%
q(z)=z^{2}+2(1+i)z$. Since $\cos (\sqrt{z})$ is an entire function, $f\in E.$
Also, $a_{0}$ and $a_{2}$ are both real and $f(q(z))=-\cosh (\pi (z+1))\in
S_{F}$. However, $f^{\prime }(0)=\dfrac{\pi }{2(1+i)}\sinh (\pi
(1+i))=\allowbreak \dfrac{\pi (i-1)}{4}\sinh \pi $, which is not real. Hence 
$f\notin S_{F}$ and $q\notin S_{F}$, but $f\circ q\in S_{F}.$

Case 2: $\pi ^{2}\notin F$

Let $f(z)=\cos (\sqrt{z+\pi ^{2}})$ and $q(z)=z^{2}+2\pi z$. Then

$f(q(z))=\allowbreak -\cos z$ $\in S_{F}$. Now $q\notin S_{F}$ since $\pi
\notin F$ and $f\notin S_{F}$ since $f^{\prime \prime }(0)=\allowbreak 
\dfrac{1}{4\pi ^{2}}\notin F$.

\QTP{Body Math}
{}


\begin{thebibliography}{9}
\bibitem{H}  Alan L. Horwitz, ``Even compositions of entire functions and
related matters'', J. Austral. Math. Soc. (Series A)63(1997), 225-237.

\bibitem{HR}  Alan L. Horwitz and Lee A. Rubel, ``When is the composition of
two power series even?'', J. Austral. Math. Soc. (Series A) 56(1994),
415-420.
\end{thebibliography}
\end{document}